\documentclass[leqno,12pt]{article}

\usepackage{amsmath}

\def\td{\mbox{\rm td}}

\newcounter{supersection}[section]
\newtheorem{pr}[supersection]{Proposition}
\newtheorem{thm}[supersection]{Theorem}

\newtheorem{re}[supersection]{Remark}

\newtheorem{ex}[supersection]{Example}

\newfont{\cirilrm}{wncyr10 scaled  1200}
\newfont{\cirilit}{wncyi10 scaled  1200}

\newcommand{\z}{\symbol{'31}}
\newcommand{\sh}{\symbol{'170}}
\newcommand{\ja}{\symbol{'37}}
\newcommand{\ijj}{\symbol{'32}}
\newcommand{\ii}{\symbol{'171}}

\def\Dj{D{\hspace{-.75em}\raisebox{.3ex}{-}\hspace{.4em}}}

\renewcommand{\thefootnote}{\arabic{footnote}}

\begin{document}

\bigskip

\centerline{\large \bf DIFFERENTIALLY TRANSCENDENTAL FUNCTIONS}

\bigskip

\bigskip

\centerline{\bf \v{Z}arko Mijajlovi\'{c}\mbox{${}^{1)}$},
Branko Male\v sevi\'{c}\mbox{${}^{2)}$}}

\bigskip

\begin{center}
{${}^{1)}$Faculty of Mathematics, University of Belgrade,           \\
$\;\,\;$Studentski trg 16, 11000 Belgrade, Serbia                   \\[2.0 ex]

${}^{2)}$Faculty of Electrical Engineering, University of Belgrade, \\
$\;\,\;$Bulevar Kralja Aleksandra 73, 11000 Belgrade, Serbia}
\end{center}

\smallskip

\begin{abstract}
The aim of this article is to exhibit a method for proving that certain analytic
functions are not solutions of algebraic differential equations.
The method is based on model-theoretic properties of differential fields and
properties of certain known transcendental differential functions,
as of $\Gamma(x)$. Furthermore,  it also determines differential
transcendence of solution of some functional equations.
\end{abstract}

\setcounter{footnote}{1}
{}\footnotetext{{\it Email address}\,:\, {\tt zarkom@eunet.yu}}
\setcounter{footnote}{2}
{}\footnotetext{{\it Email address}\,:\, {\tt malesh@eunet.yu}}

\renewcommand{\thefootnote}{}

{}\footnotetext{Second author supported in part by the project
MNTRS, Grant No. ON144020.}

\section{Notation and preliminaries}

The theory DF${}_0$ of differential fields of characteristic $0$
is the theory of fields with additional two axioms that relate to
the derivative $D$:
$$
D(x+y)= Dx + Dy,\quad D(xy)= xDy+yDx.
$$
Thus, a model of DF${}_0$ is a differential field $\mbox{\bf K} = (K,+,\cdot, D, 0,1)$
where $(K,+,\cdot,0,1)$ is a field and $D$ is a differential operator satisfying
the above axioms. A.$\,$Robinson proved that DF${}_0$ has a model completion, and then
defined DCF${}_0$ to be the model completion of DF${}_0$. Subsequently, L.$\,$Blum found
simple axioms of DFC${}_0$ without refereing to differential polynomials in more
than one variable, see \cite{Sacks72}. In the following, if not otherwise stated,
$\mbox{\bf F}, \mbox{\bf K}, \mbox{\bf L}, \ldots$ will denote differential fields,
$F,L,K,\ldots$ their domains while $\mathbf{F^\ast}, \mathbf{K^\ast}, \mathbf{L^\ast},
\ldots$ will denote their field parts, i.e. $\mathbf{F^\ast}\mbox{$=$}\mbox{$(F,+,\cdot,0,1)$}$.
Thus, $\mathbf{F^\ast}[x_1,x_2,\ldots,x_n]$ denotes the set of (ordinary) algebraic polynomials
over $\mathbf{F^\ast}$ in variables $x_1,x_2,\ldots,x_n$.  The symbol $\mbox{\bf L}\{X\}$ denotes
the ring of differential polynomials over $\mbox{\bf L}$ in the variable $X$.
Hence, if $f\in L\{X\}$ then for some $n\in N$, $N=\{0,1,2,\ldots\}$,
$f=f(X,DX,\ldots,D^nX)$ where $f\in \mathbf{F^\ast}(x,y_1,y_2,\ldots,y_n)$.

Suppose $\mbox{\bf L}\subseteq \mbox{\bf K}$. The symbol $\mbox{\rm td}(\mbox{\bf K}|\mbox{\bf L})$
denotes the transcendental degree of $\mbox{\bf K}^\ast$ over $\mathbf{L^\ast}$. The basic properties of
$\td$ are described in the following proposition.

\begin{pr}
\label{prop11}
Let $\mbox{\bf A} \subseteq \mbox{\bf B} \subset \mbox{\bf C}$
be ordinary algebraic fields. Then

\noindent
{\bf a.} $\td(\mbox{\bf B}|\mbox{\bf A})\leq
\td(\mbox{\bf C}|\mbox{\bf A})$.
\hskip 6mm
{\bf b.} $\td(\mbox{\bf C}|\mbox{\bf A})= \td(\mbox{\bf
C}|\mbox{\bf B}) + \td(\mbox{\bf B}|\mbox{\bf A})$.
\end{pr}

\noindent
The statement follows from the fact that every transcendental base
of $\mbox{\bf B}$ over $\mbox{\bf A}$ can be extended to a
transcendental base of $\mbox{\bf C}$ over $\mbox{\bf A}$, see \cite{Lang02}.

If $b\in K$, then $\mbox{\bf L}(b)$ denotes the simple differential extension of $\mbox{\bf L}$ in $\mbox{\bf K}$, i.e.
$\mbox{\bf L}(b)$ is the smallest differential subfield of $\mbox{\bf K}$ containing both $L$ and $b$. If $\mbox{\bf
L}\subseteq \mbox{\bf K}$ then $b\in K$ is {\rm differential algebraic} over $\mbox{\bf L}$ if and only if there is
a non-zero differential polynomial $f \in L\{X\}$ such that $f(b,Db,D^2b,\ldots,D^nb)=0$; if $b$ is not differential
algebraic over $\mbox{\bf L}$ then $b$ is {\rm differential transcendental} over $\mbox{\bf L}$. $\mbox{\bf K}$ is
a differential algebraic extension of $\mbox{\bf L}$ if every $b\in K$ is differential algebraic over $\mbox{\bf L}$.
Next, $f=0$ is algebraic differential equation if $f\in L\{X\}$. If $f \in L\{X\}\backslash L$, the {\rm order} of $f$,
denoted by ${\rm ord}f$, is the largest $n$ such that $D^nX$ occurs in $f$. If $f\in L$ we put $\mbox{\rm ord} \, f = -1$
and then we write $f(a) = f$ for each $a$. For $f \in L\{X\}$ we shall write occasionally $f'$ instead of $Df$,
and $f(a)$ instead of $f(a,Da,D^2a,\ldots,D^na)$ for each $a$ and $n = \mbox{\rm ord} \, f$.

Models of DCF${}_0$ are differentially closed fields.
A differential field $\mbox{\bf K}$ is {\em differentially closed} if whenever $f, g \in \mbox{\bf K}\{X\}$,
$g$ is non-zero and $\mbox{\rm ord} \, f > \mbox{\rm ord} \, g$, there is an $a \in K$ such that $f(a) = 0$ and
$g(a) \neq 0$. Let us observe that any differentially closed field is algebraically closed.
The~theory DCF${}_0$ admits elimination of quantifiers and it is submodel complete (A.$\,$Robinson):
if \mbox{$\mbox{\bf F} \subseteq \mbox{\bf L}, \mbox{\bf K}$} then \mbox{$\mbox{\bf K}_F \equiv \mbox{\bf L}_F$},
i.e.~$\mbox{\bf K}$~and~$\mbox{\bf L}$ are elementary equivalent over $\mbox{\bf F}$.

Other notations, notions and results concerning differential
fields that will be used in this article will be as in
\cite{Sacks72} or \cite{Marker96}.

\section{Differential algebraic extensions}

In this section we shall state certain properties of extensions of differential fields, which we need to
prove the main theorem 2.8. These properties, described by propositions \ref{prop21} -- \ref{prop25},
parallel in most cases properties of algebraic fields (without differential operators), by replacing
the notion of the algebraic degree with the notion of the transcendental degree of fields and ${\rm deg}f$
by ${\rm ord}f$. Proofs of these propositions  are standard and can be found in the basic literature
on differential fields and therefore they are omitted.
\begin{pr}
\label{prop21}
Suppose $\mbox{\bf L}\subseteq \mbox{\bf K}$ and let $b\in K$.
Then $b$ is differentially algebraic over $\mbox{\bf F}$ if and only if
$\td(\mbox{\bf L}(b)|\mbox{\bf L})<\infty$.
\end{pr}

\noindent
Suppose $\mbox{\bf L}\subseteq \mbox{\bf K}$ and let $a_1,a_2,\ldots,a_n\in K$ be
differentially algebraic over $\mbox{\bf L}$. Then
$$
\mbox{\bf L}(a_1,a_2,\ldots,a_i)= \mbox{\bf L}(a_1,a_2,\ldots,a_{i-1})(a_i)
$$
and $a_i$ is differentially algebraic over $\mbox{\bf L}(a_1,a_2,\ldots,a_{i-1})$ for each $i=1,2,\ldots,n$.
Thus, by propositions \ref{prop11} and \ref{prop21} we have
$$
\begin{array}{rl}
\td(\mbox{\bf L}(a_1,a_2,\ldots,a_n)|\mbox{\bf L})\,=\,
&\!\td(\mbox{\bf L}(a_1,a_2,\ldots,a_n)|\mbox{\bf L}(a_1,a_2,\ldots,a_{n-1}))\,+\hfill \\
&\!\ldots +\td(\mbox{\bf L}(a_1)|\mbox{\bf L}) < \infty,   \hfill
\end{array}
$$
so every $b\in L(a_1, \ldots,a_n)$ is differentially algebraic over $\mbox{\bf L}$.
Hence $\mbox{\bf L}(a_1,\ldots,a_n)$ is differentially algebraic over $\mbox{\bf L}$.

\begin{pr}
\label{prop22}
Let $\mbox{\bf L} \subseteq \mbox{\bf K}$ and suppose $b,a_1,\ldots,a_n\!\in\!K$. If $a_1,\ldots,a_n$
are differentially algebraic over $\mbox{\bf L}$ and $b$ is differentially algebraic over
$\mbox{\bf L}(a_1,\ldots,a_n)$, then $b$ is differentially algebraic over $\mbox{\bf L}$.
\end{pr}

\begin{pr}
\label{prop23}
If $\mbox{\bf F}\subseteq \mbox{\bf L}\subseteq \mbox{\bf K}$ and $\mbox{\bf L}$ is differentially algebraic
over $\mbox{\bf F}$ and $\mbox{\bf K}$ is differentially algebraic over $\mbox{\bf L}$ then $\mbox{\bf K}$
is differentially algebraic over $\mbox{\bf L}$.
\end{pr}

\begin{pr}
\label{prop24}
Suppose $\mbox{\bf L}\subseteq \mbox{\bf K}$ and $b\in K$. Then $b$ is a differentially algebraic over
$\mbox{\bf L}$ if and only if $Db$ is differentially algebraic over $\mbox{\bf L}$.
\end{pr}

\begin{pr}
\label{prop25}
Suppose $\mbox{\bf F}\subseteq \mbox{\bf K}$ and let
$$
L=\{b\in K\colon\, b\enskip \mbox{\it is differentially }
\mbox{\it algebraic over}\enskip \mbox{\bf F}\}.
$$

\vspace*{-2.5 mm}

\noindent
Then

\smallskip
\noindent
$\,${\bf a.}
$L$ is a differential subfield of $\mbox{\bf K}$ extending differentially algebraic $\mbox{\bf F}$.

\smallskip
\noindent
$\,${\bf b.}
If\, $\mbox{\bf K}$ is differentially algebraic closed then $\mbox{\bf L}$ is differentially algebraic closed.
\end{pr}

\begin{re}
\label{rem26}
There are differential fields with elements not having functional representation. One example of this kind is the
Hardy field of germs of functions~{\rm \cite{Gokhman96}}.
\end{re}

Let us denote by $\mathcal{M}_D$ the class of complex functions meromorphic on a
complex domain $D$. If $D=C$ then we shall write $\mathcal{M}$ instead of $\mathcal{M}_D$.
Let $\mathcal{C}$ $=$ $\mbox{\bf C}(z)$ be the differential field of complex rational functions.
Then  $\mathcal{M}$ is differential field and $\mathcal{C} \subseteq  \mathcal{M}$. Further,
let $\mathcal{L}$ $=$ $\{ f \in \mathcal{M}$ $\mbox{:} \; f \enskip \mbox{\rm differentially algebraic over}
\enskip  \mathcal{C}\}$. By Proposition \ref{prop25} $\mathcal{L}$ is a differential subfield
of $\mathcal{M}$ extending differentially algebraic $\mathcal{C}$.

\smallskip

Assume that $e$, $g$ are complex functions. If $e$ is an entire function
(i.e. holomorphic in entire finite complex plane) and $g$ is meromorphic, then
their composition $e \circ g$ is not necessary meromorphic, as the example
$e(z)=e^z$, $g(z)=1/z$ shows. However, $g \circ e$ is meromorphic since $g$ is
a quotient of entire functions, and obviously entire functions are closed under
compositions. The next proposition states that a similar property holds for
algebraic differential functions.

\begin{pr}
\label{prop27}
Let $e$, $g$ be complex algebraic differential functions over $\mathcal{C}$. If $e$ is  entire
and $g$ is  meromorphic, then  $g\circ e$ is meromorphic and differentially algebraic over $\mathcal{C}$.
\end{pr}

\noindent{\bf Proof.}
Suppose that $g$ is a solution of an algebraic differential equation over $\mathcal{C}$, i.e.
$$
f(z,g,Dg,\ldots,D^ng)\equiv 0 \;\;\mbox{for some}\,\, f \in \mathcal{C}[z,u_0,u_1,\ldots,u_n].
$$
If $e'=0$ then $e$ is a constant, so we may assume the nontrivial case, that $e'\not=0$.
There are rational expressions $\lambda_{kj}$ in $e',e'',\ldots$ such that
$$
(D^ku)\circ e= \sum_{j=1}^k \lambda_{kj}D^j(u\circ e), \quad u\,\, \text{\rm is any meromorphic function.}
$$
In fact, the sequence $\lambda_{kj}$ satisfies the recurrent identity
$$
\indent
\lambda_{k+1,j}= (\lambda_{kj}'+\lambda_{k,j-1})/e', \; \lambda_{11}=1/e',  \;
\lambda_{k0}=0, \; \lambda_{kj}=0 \; \mbox{\rm if} \; j>k.
$$

\noindent
It is easy to see that $\lambda_{kk}= (1/e')^k$. Let
$$
f_1= f(e, y, \lambda_{11}Dy,\lambda_{21}Dy+\lambda_{22}D^2y,\ldots,\lambda_{n1}Dy+\ldots+\lambda_{nn}D^ny).
$$
Then $g\circ e$ is a solution of $f_1=0$: since
$
(D^kg)\circ e= \sum_{j=1}^k \lambda_{kj}D^j(g\circ e),
$
it follows, taking $h=g\circ e$,
$$
{\begin{array}{ll}
f_1(z,h, Dh,\ldots,D^nh)&= f(e,g\circ e, (Dg)\circ e, \ldots, (D^ng)\circ e)           \cr
                        &= f(z,g,Dg,\ldots,D^ng)|_{z=e}\equiv 0.
\end{array}}
$$
Therefore, $g\circ e$ is meromorphic and differentially algebraic over $\mathcal{C}$. \hfill$\nabla$

\smallskip
Let us mention that N.$\,$Steinmetz \cite[Satz 2.]{Steinmetz80} has
considered a proposition which is the converse of the one we proved
above.

\smallskip

The next theorem is a direct corollary of the propositions \ref{prop25} and \ref{prop27}, since
$\mathcal{L}$ is a field and therefore it is closed under values of rational expressions.
\begin{thm}
\label{thm38}
Let $a(z)$ be a complex differential transcendental  function over $\mathcal{C}$,
$f(z,u_0,u_1,\mbox{\small $\ldots$},u_m,y_1,\mbox{\small $\ldots$},y_n)$ a rational expression over $\mathcal{L}$,
and assume that $e_1,\mbox{\small $\ldots$},e_m$ are entire functions which are
differentially algebraic over $\mathcal{C}$, $e_i'\not=0$, $1\leq i\leq m$. If $b$ is meromorphic
and $f(z,b,b\circ e_1,\ldots,b\circ e_m,Db,\ldots,D^nb)\equiv a(z)$ then $b$
is differential transcendental over $\mathcal{C}$.
\end{thm}

\medskip
We shall use it in proofs of differential transcendentality of
certain complex functions. In most cases functions $e_i(z)$ will be linear
functions $\alpha z+ \beta$, $\alpha\not=0$, and $f$ will be a polynomial over
$\mathcal{C}$ (observe that $\mathcal{C} \!\subset\! \mathcal{L}$).
Also, it is possible to consider differential transcendental functions
when $f$ is a rational expression over $\mathcal{L}$, for example as a solution
of the equation \mbox{$y''\!(z) + y'\!(z)/y(z) + y(\sin z) = \Gamma(z\!+\!1)$}.

\section{Differentially transcendental functions}

Suppose $\mbox{\bf L}\subseteq \mbox{\bf K}$. Let $\mbox{\bf R}(x)$ be the differential field
of real rational functions. The following H\"older's famous theorem asserts the differential
transcendentality of Gamma function, see \cite{Holder1887}.

\begin{thm}
\label{thm31}
{\bf a.} $\Gamma(x)$ is not differentially algebraic over $\mbox{\bf R}(x)$.\enskip
{\bf b.} $\Gamma(z)$ is not differentially algebraic over $\mbox{\bf C}(z)$.~$\nabla$
\end{thm}

Now we shall use the differential transcendence of $\Gamma(z)$ and
properties of differential fields developed in the previous  section
to prove differential transcendentality over $\mathcal{C}$  of some
analytic functions. $\Gamma(z)$ is  meromorphic and by H\"older's theorem
and Proposition \ref{prop24}, $\Gamma(z)\not\in\mathcal{L}$.

\noindent
\begin{ex}
\label{exmp32}
{\rm This example is an archetype of proofs,
based on properties of differential fields exhibited in the previous section,
of differential transcendentality of some well known complex functions.
The Riemann zeta function is differentially transcendental over $\mathcal{C}$ (Hilbert).
First we observe that $\zeta(s)$ is  meromorphic  and that $\zeta(s)$
satisfies the well-known functional equation:
$$
\zeta(s) = \chi(s)\zeta(1-s),\quad \mbox{\rm where}\quad \chi(s) =
\displaystyle\frac{(2\pi)^s}{
2 \Gamma(s) \cos (\mbox{\small $\displaystyle\frac{\pi s}{2}$}) }.
$$
Now, suppose that $\zeta(s)$ is differentially algebraic over $\mathcal{C}$, i.e. that
$\zeta(s)\in \mathcal{L}$. Then $\zeta(1-s)$ and $\zeta(s)/\zeta(1-s)$ belong to $\mathcal{L}$,
too, so $\chi(s)$ belongs to $\mathcal{L}$. The elementary functions $(2\pi)^s$, and
$\cos({\frac{\pi s}{2}})$ obviously are differentially algebraic over $\mathcal{C}$
i.e. they belong to $\mathcal{L}$. As $\mathcal{L}$ is a field, it follows that
$\Gamma(z)$ is differentially algebraic over $\mathcal{C}$. Hence,
$\Gamma(z)\in \mathcal{L}$ but this  yields a contradiction.
Therefore $\zeta(s)$ is differentially transcendental function over $\mathcal{C}$.
Generally, Dirichlet $L$-series
$$
L_{k}(s)= \displaystyle\sum\limits_{n=1}^{\infty}{\chi_{k}(n)
\displaystyle\frac{1}{n^s}} \quad (k \in \mbox{\bf Z}),
$$
where $\chi_{k}(n)$ is Dirichlet character see \cite{IrelandRosen82},
is differentially transcendental function over $\mathcal{C}$.
This follows from a well-known functional equations
$$
L_{-k}(s) = 2^s \pi^{s-1} k^{-s+\frac{1}{2}} \Gamma(1-s) \cos (\mbox{\small
$\displaystyle\frac{\pi s}{2}$}) L_{-k}(1-s)
$$
and
$$
L_{+k}(s) = 2^s \pi^{s-1} k^{-s+\frac{1}{2}} \Gamma(1-s) \sin (\mbox{\small
$\displaystyle\frac{\pi s}{2}$}) L_{+k}(1-s),
$$
where $k \in \mbox{\bf N}$. D.$\,$Gokhman also proved in \cite{Gokhman96},
but in entirely different way, that Dirichlet series are differentially
transcendental. Besides Riemann zeta function $\zeta(s)=L_{+1}(s)$,
Dirichlet eta function
$$
\eta(s) =
\displaystyle\sum\limits_{n=1}^{\infty}{(-1)^{n+1}\displaystyle\frac{1}{n^s}}
= (1-2^{1-s})L_{+1}(s)
$$
and Dirichlet beta function
$$
\beta(s) =
\displaystyle\sum\limits_{n=0}^{\infty}{(-1)^{n}\displaystyle\frac{1}{(2n+1)^s}}
= L_{-4}(s)
$$
are transcendental differential functions as examples of Dirichlet
series see \cite[page 289]{Borwein87}. \hfill$\nabla$}
\end{ex}
\begin{ex}
\label{exmp33}
{\rm Kurepa's Function related to the Kurepa Left Factorial Hypothesis,
see \cite[problem B44]{Guy94}, is defined by \cite{Kurepa73}:
\begin{equation}
\label{K_INT_1}
K(z) = \displaystyle\int\limits_{0}^{\infty}{
e^{-t} \displaystyle\frac{t^{z}-1}{t-1} \: dt}.
\end{equation}
Kurepa established in \cite{Kurepa73} that $K(z)$ can be continued
meromorphically to the whole complex plane. Also, it satisfies the
recurrence relation
\begin{equation}
\label{K_FE_1}
K(z) - K(z-1) = \Gamma(z),
\end{equation}
hence, by Theorem \ref{thm38},  $K(z)$ is differential transcendental over $\mathcal{C}$.
The functional equation (\ref{K_FE_1}), besides Kurepa's function $K(z)$, has another
solution which is given by the  series
\begin{equation}
\label{Def_K1}
K_{1}(z) =
\displaystyle\sum\limits_{n=0}^{\infty}{\Gamma(z-n)}.
\end{equation}
This function has simple poles at integer points. We can conclude similarly that $K_1(z)$
is transcendental differential function. Between functions $K(z)$ and $K_1(z)$ the
following relation is true see \cite{Slavic73} and \cite{Malesevic03}:
$$
K(z) = \displaystyle\frac{\mbox{\rm Ei}(1)}{e} -
\displaystyle\frac{\pi}{e} \, \mbox{\rm ctg}{\, \pi z} + K_{1}(z),
$$
where $\mbox{\rm Ei}(x)$ is the exponential integral see \cite{GradsteinRyzhik71}.
Each meromorphic solution of the functional equation {\rm (\ref{K_FE_1})}
is differential transcendental over $\mathcal{C}$.

Alternating Kurepa's function related to the alternating sums  of factorials
see \cite[problem B43]{Guy94}, is defined by  \cite{Petojevic02}:
\begin{equation}
\label{A_INT_1}
A(z) = \displaystyle\int\limits_{0}^{\infty}{
e^{-t} \displaystyle\frac{t^{z+1}-(-1)^{z}t}{t+1} \: dt}.
\end{equation}
$A(z)$ is meromorphic  and it satisfies the functional equation
\begin{equation}
\label{A_FE_1}
A(z) + A(z-1) = \Gamma(z+1).
\end{equation}
The functional equation (\ref{A_FE_1}), besides alternating Kurepa's function $A(z)$,
has another solution which is given by the  series
\begin{equation}
\label{Def_A1}
A_{1}(z) =
\displaystyle\sum\limits_{n=0}^{\infty}{(-1)^n\Gamma(z+1-n)}.
\end{equation}
As in the case of $K(z)$, we see that $A(z)$ and $A_1(z)$ are differential transcendental over
$\mathcal{C}$. The following identity holds, see \cite{Malesevic04}:
$$
A(z) = - {\big (}1 + e \, \mbox{\rm Ei}(-1){\big )}(-1)^z +
\displaystyle\frac{\pi e}{\sin \pi z} + A_{1}(z).
$$
Generally, each meromorphic solution of a functional equation {\rm
(\ref{A_FE_1})} is transcendental differential function over the
field $\mathcal{C}$.

G.$\,$Milovanovi\'{c} introduced in \cite{Milovanovic96} a sequence
of meromorphic functions
$$
K_{m}(z) - K_{m}(z-1) = K_{m-1}(z), \quad K_{-1}(z) = \Gamma(z),
\; \; K_{0}(z) = K(z).
$$
By use of induction and Theorem \ref{thm38} we can conclude that  $K_m(z)$ are differential transcendental
over $\mathcal{C}$. Analogously, a sequence of meromorphic functions defined by
$$
A_{m}(z) + A_{m}(z-1) = A_{m-1}(z), \quad A_{-1}(z) = \Gamma(z+1),
\; \; A_{0}(z) = A(z)
$$
is a sequences of transcendental differential functions over
$\mathcal{C}$. \hfill$\nabla$}
\end{ex}

\noindent
\begin{ex}
\label{exmp34}
{\rm The meromorphic function
$$
 H_{1}(z)= \sum_{n=0}^{\infty}{\displaystyle\frac{1}{(n+z)^2}}
$$
is differentially transcendental over $\mathcal{C}$. Really, $D^2\ln(\Gamma(z))=H_{1}(z)$, i.e.
$\Gamma(z)$ satisfies the algebraic differential equation $(D^2\Gamma)\Gamma-(D\Gamma)^2-H_{1}\Gamma^2=0$
over $\mathcal{C}(H_{1})$. Thus, if $H_{1}$ would be differentially algebraic over $\mathcal{C}$,
then by Proposition \ref{prop23}, $\Gamma$ would be too, what yields a contradiction. Hence,
$H_{1}(z)$ is differentially transcendental function over $\mathcal{C}$.  \hfill$\nabla$}
\end{ex}
\begin{re}
\label{rem35}
{\rm We see that in Example \ref{exmp34} functions $\Gamma(z)$ and
$H_{1}(z)$ are differentially algebraically dependent, i.e. $g(\Gamma,H_{1})=0$,
where $g(x,y)= x''x-(x')^2-yx^2$. We do not know if similar
dependencies exist for pairs $(K,\Gamma)$ and $(\zeta,\Gamma)$.
It is very likely that these pairs are in fact differentially
transcendental.}
\end{re}
\begin{ex}
\label{exmp36}
{\rm W.$\,$Gautschi and J.$\,$Waldvogel considered in \cite{Gautschi01}
a family of functions which satisfies
$$
I_\alpha(x) = x I_{\alpha-1}(x) + \Gamma(\alpha), \quad \alpha > -1, \;x > 0.
$$
For each $x_0 > 0$ function $f(\alpha) = I_{\alpha}(x_0)$ can be
continued meromorphically  (with respect to the complex $\alpha$,
see \cite{Gautschi01}).  By Theorem \ref{thm38} the function $f(\alpha)$
is differential transcendental over $\mathcal{C}$. \hfill$\nabla$}
\end{ex}
\begin{ex}
\label{exmp37}
{\rm The Barnes $G$-function is defined by \cite{Barnes1900}:
$$
G(z+1) = (2\pi)^{z/2}e^{-(z(z+1) + \gamma z^2)/2}
\displaystyle\prod\limits_{n=1}^{\infty}{ {\bigg (} {\Big(} 1 +
\displaystyle\frac{z}{n} {\Big )}^{n} e^{-z+z^2/(2n)} {\bigg )}},
$$
where $\gamma$ is the Euler constant. Barnes $G$-function is an
entire function and it satisfies following functional equation
\begin{equation}
\label{G_FE_1}
G(z+1) = \Gamma(z) G(z),
\end{equation}
By Theorem \ref{thm38}  each meromorphic solution of the functional equation (\ref{G_FE_1})
is differential transcendental over $\mathcal{C}$. \hfill$\nabla$}
\end{ex}
\begin{ex}
\label{exmp38}
{\rm The Hadamard's factorial function is defined by:
$$
H(z) = \displaystyle\frac{1}{\Gamma(1-z)}
\displaystyle\frac{d}{dz} \ln {\bigg (} \Gamma {\Big (}
\displaystyle\frac{1-z}{2} {\Big )} \mbox{\Large $/$} \Gamma {\Big
(} 1-\displaystyle\frac{z}{2} {\Big )} {\bigg )},
$$
see \cite{Davis59}. Hadamard's function is an entire function and
it satisfies following functional equation
\begin{equation}
\label{H_FE_1}
H(z+1) = zH(z) +
\displaystyle\frac{1}{\Gamma(1-z)}.
\end{equation}
By Theorem \ref{thm38} each meromorphic solution of  the functional equation (\ref{H_FE_1})
is differential transcendental over $\mathcal{C}$. \hfill$\nabla$}
\end{ex}
\begin{ex}
\label{exmp39}
{\rm Let $\tau(n)$ be defined by $\displaystyle\sum\limits_{n=1}^{\infty}{\tau(n)q^n} = q
\displaystyle\prod\limits_{n=1}^{\infty}{(1-q^n)^{24}}, \; q=e^{2 \pi i z}.$ This function
was first studied by \cite{Ramanujan16}. Ramanujan's Dirichlet $L$-series defined by
$$
 f(s) = \displaystyle\sum\limits_{n=1}^{\infty}{\displaystyle\frac{\tau(n)}{n^s}},
$$
is differentially transcendental function over $\mathcal{C}$. This follows from functional
equation see \cite[page 173]{Hardy59}:

\smallskip

\centerline{$
f(s) = \displaystyle\frac{(2 \pi)^s}{(2\pi)^{12-s}} \,
\displaystyle\frac{\Gamma(12-s)}{\Gamma(s)} \, f(12-s)
$}

\vspace*{-2.5 mm}

\noindent
i.e.

\centerline{$
\qquad\qquad
f(s)
=
-\displaystyle\frac{(2 \pi)^s}{(2\pi)^{12-s}}
\, \displaystyle\frac{\pi}{\sin ( \pi s )}
{\bigg (} \displaystyle\prod_{i=1}^{11} {(s-i)} {\bigg )}
\, \displaystyle\frac{1}{\Gamma(s)^2} \, f(12-s).
\qquad\qquad\nabla$}}
\end{ex}
\begin{ex}
\label{exmp310}
{\rm Let us note that function \cite{GradsteinRyzhik71} (integral 3.411/1):

\smallskip

\centerline{$\displaystyle\int\limits_{0}^{\infty}{
\displaystyle\frac{x^{p-1}}{e^x-1} \, dx} = \Gamma(p) \zeta(p)$}

\noindent
is differentially transcendental over $\mathcal{C}$. This follows from the functional equation
\mbox{$\zeta(1\!-\!s)$} \mbox{$=2(2\pi)^{-s}\cos(\frac{s \pi}{2}) {\big (}\Gamma(s)\zeta(s){\big )}$}.
Analogously function $\Gamma(p) \beta(p)$ is differentially transcendental over $\mathcal{C}$.
On the basis of differentially transcendental functions $\Gamma(p), \Gamma(p) \zeta(p), \Gamma(p) \beta(p)$
it is possible to find large number of examples differentially transcendental functions
from table of integrals \cite{GradsteinRyzhik71} (e.g. integrals: 3.411/3, 3.523/1,
3.523/3, 3.551/2, 3.769/1, 3.769/2, 3.944/5, 3.944/6). \hfill$\nabla$}
\end{ex}

\end{document}